\documentclass[a4paper,11pt]{article}
\usepackage{makeidx}

\usepackage[T1]{fontenc}
\usepackage[dvips]{graphicx}
\usepackage{amsmath,amscd,amsthm}
\usepackage{amssymb}
\usepackage{tabularx}
\usepackage{graphics}
\usepackage[latin1]{inputenc}
\graphicspath{c:/MISDOC~1}\makeindex \textwidth 15cm \textheight
22cm \topmargin 0.5cm
\usepackage[dvipdfm,dvipsnames,usenames]{color}
\usepackage{soul}

\newtheorem{prp}{Proposition}

\newcommand{\N}{\mathbb{N}}

\newcommand{\C}{\mathbb{C}}

\begin{document}

\begin{center}
{\Large \bf  Minimal Faithful Upper-Triangular Matrix
Representations for Solvable Lie Algebras}
\end{center}

\begin{center}
{\bf Manuel Ceballos$^\dagger$, Juan N\'u\~nez$^\dagger$ and \'Angel F. Tenorio$^\ddagger$\footnote{Corresponding author. Phone: +34-954349354. Fax: +34-954349339}}
\end{center}

\begin{center}
\small $^\dagger$Departamento de Geometr\'{\i}a y Topolog\'{\i}a.\\
Facultad de Matem\'aticas. Universidad de Sevilla.\\
Aptdo. 1160. 41080-Seville (Spain).\\
Email: \{mceballos,jnvaldes\}@us.es\\
\quad \\
$^\ddagger$Dpto. de Econom\'{\i}a, M\'etodos Cuantitativos e H.$^{\rm a}$ Econ\'omica.\\
Escuela Polit\'ecnica Superior. Universidad Pablo de Olavide.\\
Ctra. Utrera km. 1, 41013-Seville (Spain).\\
Email: aftenorio@upo.es
\end{center}
\vspace{0.25cm}

\begin{quote}
\noindent {\bf Abstract.} A well-known result on Lie Theory states that every finite-dime\-nsional complex solvable Lie algebra can be represented as a matrix Lie algebra, with upper-triangular square matrices as elements. However, this result does not specify which is the minimal order of the matrices involved in such representations. Hence, the main goal of this paper is to revisit and implement a method to compute both that minimal order and a matrix representative for a given solvable Lie algebra. As application of this procedure, we compute representatives for each solvable Lie algebra with dimension less than $6$.
\end{quote}

\vspace{0.2cm}



\section{Introduction}

Representation Theory of Lie algebras can be allowed for the classification of Lie algebras and groups, which has broad applications to the analysis of continuous symmetries in Mathematics and Physics. More concretely, in Mathematics, the classification of Lie groups reveals symmetries in differential equations. With respect to Physics, representation theory yields natural connections between representation of Lie algebras and the properties of elementary particles.

Ado's Theorem states that given a finite-dimensional complex Lie algebra $\mathfrak{g}$, there exists a matrix algebra isomorphic to $\mathfrak{g}$ (see \cite{Jac} for the classical proof and  \cite{Ner} for a very short alternative). In this way, every finite-dimensional complex Lie algebra can be represented as a Lie subalgebra of the complex general linear algebra $\mathfrak{gl}(n;{\C})$, of complex $n \times n$ matrices, for some $n \in {\N}$.

This paper focuses on Lie algebra $\mathfrak{h}_n$, of $n\times n$ upper-triangular matrices. It is well known that every finite-dimensional solvable Lie algebra is isomorphic to a subalgebra of $\mathfrak{h}_n$, for some $n\in {\N}$ (see \cite[Proposition 3.7.3]{Var98}). Therefore, the following interesting question arises for a given finite-dimensional solvable Lie algebra $\mathfrak{g}$: determining the minimal $n \in \mathbb N$ such that $\mathfrak h_n$  contains $\mathfrak g$ as a Lie subalgebra; i.e. obtaining the minimal faithful representation of $\mathfrak g$ by using $n\times n$ upper-triangular matrices.

Several authors have studied the minimal dimension $\mu(\mathfrak g)$ to represent a given Lie algebra $\mathfrak g$ (see Burde \cite{Bu98}
, for instance). However, most of them have considered faithful $\mathfrak g$-modules instead of the particular subclass consisting of Lie algebras~$\mathfrak h_n$. Therefore, the value of $\mu(\mathfrak g)$ is less than or equal to the dimension to be computed in this paper. Regarding this matter, matrix representations were computed by Ghanam et al. \cite{G} for low-dimensional nilpotent Lie algebras, but not studying the minimality and giving some non-minimal representations.

The interest in these faithful representations is motivated, among other issues, by problems from Geometry and Topology. For example, Milnor \cite{Mil77} and Auslander \cite{Aus64,Aus77} studied generalizations of crystallographic groups in relation with this minimal value for matrix representations. Another motivation is based on the following result: Lie algebra $\mathfrak{g}$ of a given Lie group $G$ admitting a left-invariant affine structure satisfies that its minimal dimension of faithful representations is $\mu(\mathfrak{g}) \leq n+1$.

Several papers throughout the literature deal with matrix representation of these solvable Lie algebras. For example, Benjumea et al. \cite{B} introduced an algorithmic procedure which explicitly computed a representative of minimal faithful unitriangular matrix representations for a given nilpotent Lie algebra and its associated Lie group. Subsequently, the complete list of minimal faithful unitriangular matrix representations was given by Benjumea et al.~\cite{MathScan08} for nilpotent Lie algebras of dimension less than 6. Finally, the matrix representation of filiform Lie algebras of dimension less than $9$ was computed in \cite{CNT}. Additionally, N\'u\~nez and Tenorio \cite{NuTeMega07} continued with this research and gave the outlines of an algorithmic procedure to compute explictly representatives of the minimal faithful matrix representation for solvable Lie algebras by using Lie algebras $\mathfrak h_n$, giving some examples of application by hand. This procedure adapted that given in Benjumea et al. \cite{B}, but neither the algorithm was completely debugged nor implementations were carried out and run.

The main goal of the current paper is to advance in the above-mentioned research by debugging and implementing the algorithm sketched in \cite{NuTeMega07} in order to automate the computation of minimal faithful matrix representations for a given solvable Lie algebra starting from its law. As application, we have also computed representations for each solvable Lie algebra of dimension less than $6$ as well as for others of higher dimension. To do so, we have used the classifications given by Mubarakzyanov and Turkowski (see \cite{AND,Mub,Mub1,Turi}).

This paper is structured as follows: Section $2$ reviews some well-known results on Lie Theory to be applied later. Thereupon, Section $3$ revisits the algorithmic procedure sketched in \cite{NuTeMega07} to compute minimal faithful representations for solvable Lie algebras by using upper-triangular matrices, incorporating a formulation of the algorithm which can be dealt computationally with implementation in MAPLE 12. To shorten the paper length, the computational method is only explicitly applied to two algebras in Section $4$. Just afterwards, Section $5$ gives an explicit list with representatives of minimal faithful matrix representations for solvable Lie algebras of dimension less than $6$.

\section{Preliminaries} \label{Sec1}

For an overall review  on Lie algebras,  the reader can consult \cite{Var98}. In the present section, we only recall some definitions and results about Invariant Theory and Lie algebras to be applied later. Throughout this article, we only consider finite-dimensional Lie algebras over the complex number field $\mathbb C$.

Given a Lie algebra $\mathfrak g$, its {\it derived series} is defined as follows
\begin{equation}\mathcal{C}_1(\mathfrak g)= \mathfrak g, \ \mathcal{C}_2(\mathfrak g)=[\mathfrak g,\mathfrak g], \ \mathcal{C}_3(\mathfrak g)=[\mathcal{C}_2(\mathfrak g),\mathcal{C}_2(\mathfrak g)], \ \dots, \ \mathcal{C}_k(\mathfrak g)=[\mathcal{C}_{k-1}(\mathfrak g), \mathcal{C}_{k-1}(\mathfrak g)], \ \dots \end{equation}
Additionally, the Lie algebra $\mathfrak g$ is said to be {\it solvable} if there exists a natural integer $m$ such that $\mathcal{C}_m(\mathfrak g) \equiv 0$. The solv-index of $\mathfrak g$ is precisely the value of $m\in\mathbb N$ such that $\mathcal{C}_m(\mathfrak g) = 0$ and $\mathcal{C}_{m-1}(\mathfrak g) \neq 0$.

The relation between the derived series of a given Lie algebra $\mathfrak g$ and that of a Lie subalgebras is given as follows

\begin{prp}\label{PropCon}
If $\mathfrak h$ is a Lie subalgebra of a given Lie algebra $\mathfrak g$, then $\mathcal{C}_{k}(\mathfrak h)\subseteq
\mathcal{C}_{k}(\mathfrak g)$, for all $k \in \mathbb N$.
\end{prp}

Given $n\in {\N}$, the complex solvable Lie algebra $\mathfrak{h}_n$ consists of $n \times n$ upper-triangular matrices; i.e. its vectors are expressed as
\begin{equation}
h_{n}(x_{r,s})=\left(
\begin{array}{cccc}
x_{11} & x_{12} & \cdots & x_{1n} \\
0 & x_{22} & \cdots & x_{2n} \\
\vdots & \vdots & \ddots & \vdots \\
0 & \cdots & 0 & x_{nn}
\end{array} \right), \qquad \mathrm{with \ } x_{r,s} \in \mathbb C, \ \mathrm{for} \ 1\le r\le s\le n.
\end{equation}

Lie algebra $\mathfrak{h}_n$ has a basis $\mathcal B_n$ consisting of vectors $X_{i,j}=h_n(x_{r,s})$ with $1 \leq i \leq j \leq n$ and such
that
\begin{equation}
x_{r,s}=\left\{\begin{array}{ll} 1, & \mbox{if
}(r,s) = (i,j),\\ 0, & \mbox{if }(r,s)\neq (i,j).
\end{array} \right.
\end{equation}
\vspace{-0.7em}

The dimension of $\mathfrak{h}_n$ is $\frac{n(n+1)}{2}$ and the nonzero brackets with respect to basis $\mathcal B_n$ are
\begin{eqnarray}
\label{EcBrack1} [X_{i,j},X_{j,k}] = X_{i,k}, & \ \forall\ 1 \leq i < j < k \leq n;\\
\label{EcBrack2} {[}X_{i,i}, X_{i,j}] = X_{i,j}, & \ 1 \leq i < j \leq n;\\
\label{EcBrack3} {[}X_{k,i}, X_{i,i}] = X_{k,i}, & \ \forall\ k \leq i \leq n.
\end{eqnarray}\vspace{-0.5em}

\section{Computing Minimal Matrix Representations}
\label{Algorithmic method}

This section continues the work started in \cite{NuTeMega07} and introduces an algorithmic method to compute minimal matrix representations of solvable Lie algebras in such a way this can be dealt with computer algebra. After explaining step by step the algorithm, this is implemented in Maple 12 and applied to several examples.

Given a Lie algebra $\mathfrak g$, a {\it representation} of $\mathfrak g$ in $\mathbb C^n$ is a homomorphism of Lie algebras  $\phi:\mathfrak g \rightarrow \mathfrak{gl}(\mathbb C,n)$. Then $n\in\mathbb N$ is called the {\it dimension} of this representation. Ado's theorem states that every finite-dimensional Lie algebra over a field of characteristic zero  (as in the case of $\mathbb C$) has a linear injective representation on a finite-dimensional vector space; that is, a {\it faithful representation}.

Usually, representations are defined as {\it $\mathfrak g$-modules}, consisting of homomorphisms of Lie algebras from $\mathfrak g$ to Lie algebra $\mathfrak{gl}(V)$ of endomorphisms over an arbitrary $n$-dimensional vector space $V$ (like in \cite{F-H91}).

Regarding minimal representations of Lie algebras, Burde \cite{Bu98} introduced the invariant $\mu(\mathfrak g)$ for an
arbitrary Lie algebra $\mathfrak g$
$$\mu(\mathfrak g)={\rm min}\{{\rm dim} (M) \ | \ M \mbox{ is a faithful } \mathfrak g\mbox{-module}\}.$$
In this section, matrix faithful representations of solvable Lie algebras are studied. Moreover, we are interested in minimal faithful matrix representations with a particular restriction: the representation is contained in $\mathfrak h_n$ for some $n\in\mathbb N$. In this way, given a solvable Lie algebra $\mathfrak g$, we want to compute the minimal value $n$ such that $\mathfrak h_n$ contains  a Lie subalgebra isomorphic to $\mathfrak g$.  This value is  also  an invariant of $\mathfrak g$ and its expression is given by
$$\bar{\mu}(\mathfrak g)={\rm min}\{n\in \mathbb N \ | \ \exists \mbox{ subalgebra of } \mathfrak h_n \, \mbox{ isomorphic to } \mathfrak g\}.$$
In general, invariants $\mu(\mathfrak g)$ and $\bar{\mu}(\mathfrak g)$ can be different from each other.

Next, we show the algorithmic method to compute minimal faithful matrix representations for those algebras by using Lie algebras $\mathfrak h_n$. The minimality must be understood in the following sense: There exists a faithful matrix representation of $\mathfrak g$ in $\mathfrak{h}_n$, but no in $\mathfrak{h}_{n-1}$.

To do  so, we give a  step-by-step  explanation of the algorithm used to determine these minimal representations for a given solvable Lie
algebra $\mathfrak{g}$ of dimension $n$.

\begin{enumerate}

\item According to Proposition \ref{PropCon}, we compute the first natural integer $k$ such that the derived series of $\mathfrak{h}_k$ is compatible with that associated with~$\mathfrak g$.

\item  We search  a Lie subalgebra of $\mathfrak h_k$ isomorphic to $\mathfrak g$, with $k$ as low as possible. To do so, the vectors in the basis $\{e_i\}_{i=1}^n$  of $\mathfrak g$ are expressed as the following linear combinations of basis $\mathcal{B}_k$
    \begin{equation}\label{comblineal} e_h = \sum_{1\le i \le j\le k} \lambda_{i,j}^h X_{i,j}, \quad {\rm for} \,\, 1 \leq h \leq n. \end{equation}

\item Bracket $[e_i,e_j]$ is computed for $1\le i\le j\le n$. When imposing the law of $\mathfrak{g}$, a system of non-linear equations is obtained by comparing coordinate to coordinate with respect to basis $\mathcal B_k$.

\item We solve the system of equations  
    and a solution of the system provides us one of the representations searched for Lie algebra $\mathfrak{g}$ if the solution corresponds to a set of vectors being linearly independent. When no solution is obtained, Lie algebra $\mathfrak{g}$ cannot be represented as a Lie subalgebra of $\mathfrak{h}_k$. In this case, we go back to Step 2 and repeat each step with Lie algebra $\mathfrak{h}_{k+1}$.
\end{enumerate}

The representation obtained for Lie algebra $\mathfrak{g}$ is minimal because we start with $k=1$ and $k$ increases one unit when no representation can be obtained from $\mathfrak{h}_k$.

Obviously, the set of solutions of the system in the last step depends on the Lie algebra which the algorithm is applied to. As an example of application, we have minimally represented solvable Lie algebras with dimension less than $6$ in Section \ref{Representation}. The set of solutions in Step~4 has been computed with the command {\tt solve} in the symbolic computation package Maple $12$. This command works efficiently with polynomial equations, receives as inputs the list of equations and the list of variables, and returns as output the algebraic expression of the set of solutions.


Furthermore, in order to compute a particular solution of the previous system, we have searched one having as many coefficients $\lambda_{i,j}^h$ being equal to 0 as possible. In this way, a coefficient is assumed to be equal zero when it does not appear in the relations obtained by the equations.  This will be a natural representative of the Lie algebra $\mathfrak g$.

\subsection{Implementation}

Next, we show the implementation of the different routines in order to apply the previous method. They have been written using the symbolic computation package MAPLE $12$, loading the libraries {\tt DifferentialGeometry, LieAlgebras} to activate commands related to Lie algebras.

First, the routine {\tt law$\_$h} is implemented to compute the law of the solvable Lie algebra $\mathfrak{h}_n$. This routine receives as input the value of $n$ and returns the list of brackets expressing the law of $\mathfrak{h}_n$ with respect to the basis $\{e_1,e_2, \ldots, e_{\frac{n (n+1)}{2}}\}$, which corresponds to $\{X_{1,1}, \ldots, X_{1,n}, X_{2,2}, \ldots, X_{2,n}, \ldots, X_{n,n}\}$. For the implementation, a list {\tt B} saves the basis of $\mathfrak{h}_n$ and {\tt S} keeps all the (non-zero) brackets involved in the law. To carry out the computations, three different loops are programmed to find and save the three different types of non-zero brackets in Eqs.~(\ref{EcBrack1})--(\ref{EcBrack3}). Finally, the law of the algebra is saved in the variable {\tt Ext1} to be  loaded  in a later routine.

{\footnotesize
\begin{verbatim}
> law_h:=proc(n)
> local B, S;
> B:=[]; S:=[];
> for i from 1 to n do                                   (*Constructing the basis*)
>   for j from i to n do
>     B:=[op(B),X[i,j]];
>   end do;
> end do;
> for i from 1 to n-1 do                                 (*Finding brackets in Eq. (5)*)
>   for j from i+1 to n do
>     S:=[op(S),[X[i,i],X[i,j]]=X[i,j]];
>   end do;
> end do;
> for i from 1 to n-1 do                                 (*Finding brackets in Eq. (6)*)
>   for j from i+1 to n do
>     S:=[op(S),[X[i,j],X[j,j]]=X[i,j]];
>   end do;
> end do;
> for i from 1 to n-2 do                                 (*Finding brackets in Eq. (4)*)
>   for j from i+1 to n-1 do
>     for k from j+1 to n do
>       S:=[op(S),[X[i,j],X[j,k]]=X[i,k]];
>     end do;
>   end do;
> end do;
> return LieAlgebraData(S,B,Ext1, "LieAlgebraData");     (*Defining the algebra*)
> end proc:
\end{verbatim}}

Next, the routine {\tt DerviedSeries$\_$h} receives as input the value of $n$ and computes a list with the dimension of each term in the derived series of $\mathfrak{h}_n$. Let us note that we have to distinguish two different cases; being the first when $n$ is very low. Otherwise, we only have to consider powers of 2 as it can be proved by a straightforward inductive reasoning.

{\footnotesize
\begin{verbatim}
> DerivedSeries_h:=proc(n)
> local L;
> L:=[n*(n+1)/2];
> if n<4 then
>   for i from 1 to n do
>     L:=[op(L),(n-i)*(n-i+1)/2];
>   end do;
> else
    for i from  0 by 1 while 2^i<n do
>     L:=[op(L),(n-2^i)*(n-2^i +1)/2];
>   end do;
> end if;
> if member(0,L)=false then L:=[op(L),0]; end if;
> return L;
> end proc:
\end{verbatim}}

%
%

Now, we define solvable Lie algebra $\mathfrak{g}$ according to the following notation
{\footnotesize
\begin{verbatim}
> L:= _DG([["LieAlgebra", g, [n]], [A]]);
> DGsetup(L);
\end{verbatim}
}
\noindent where $n$ is the dimension of $\mathfrak{g}$ (a value inserted by  the user) and {\tt A} is a list containing information about the structure constants of the law of $\mathfrak{g}$. Elements in {\tt A} must be of the form {\tt [[i, j, k], cijk]} where {\tt cijk} is the structure constant $\lambda_{i,j}^k$ corresponding to the coefficient of $e_k$  in the bracket $[e_i,e_j]$. Once these data are loaded, we can operate over Lie algebra $\mathfrak{g}$.

The next routine is called {\tt DerivedSeries} and computes a list with the dimension of each ideal in the derived series of $\mathfrak{g}$.

{\footnotesize
\begin{verbatim}
g > DerivedSeries:=proc()
g > local k;
g > C[1]:=[seq(e||i, i=1..n)]; C[2]:=DerivedAlgebra();        (*Initiating derived series*)
g > if C[2]=[] then return "Abelian Lie algebra"; end if;     (*Testing abelian Lie algebras*)
g > for i from 3 by 1 while C[i]<> [] do                      (*Constructing derived series*)
g >   C[i]:=DerivedAlgebra(C[i-1]);
g >   if C[i]=[] then return [seq(nops(C[j]),j=1..i)]; end if;
g > end do;
g > end proc:
\end{verbatim}
}

Next, the solv-index of $\mathfrak{g}$ is saved in the variable called {\tt index$\_$g} and all the ideals of the derived series are also defined. To do so, we execute the following sentences


{\footnotesize
\begin{verbatim}
g > assign(index_g,nops(DerivedSeries()));
g > C[1](g):=[seq(e||i, i=1..n)];
g > C[2](g):=DerivedAlgebra();
g > for i from 3 to index_g do
g >   C[i](g):=DerivedAlgebra(C[i-1](g));
g > end do;
\end{verbatim}}

%

The routine {\tt DimRepresentation} computes the minimal dimension for a matrix representation of $\mathfrak{g}$ by using Lie algebras $\mathfrak{h}_n$. To implement this routine, we compare the dimension sequence of the derived series of both $\mathfrak{g}$ and $\mathfrak{h}_n$ by using the routines {\tt DerivedSeries$\_$h} and {\tt DerivedSeries}. The output is the minimal {\tt k}$\in\mathbb N$ such that Step 1 is verified.

{\footnotesize
\begin{verbatim}
g > DimRepresentation:=proc()
g > L:=DerivedSeries();
g > k:=0; d:=nops(L);
g > if d <= 4 then                 (*Initiating dimension of representation*)
g >   n:=d-1; else n:=d;
g > end if;
g > while k=0 do
g >   M:=DerivedSeries_h(n);
g >   for i from 1 to nops(L) do
g >     if L[i]<=M[i] then k:=k;   (*Comparing dimension sequences of derived algebras*)
g >       else k:=k+1;
g >     end if;
g >   end do;
g >   if k>0 then k:=0;            (*Checking compatibility between dimension sequences*)
g >     else return n;             (*Returning dimension of representation*)
g >   end if;
g >   n:=n+1;
g > end do;
g > end proc:
\end{verbatim}}

After completing Step 1 in the algorithm, we need to express all the vectors in the basis of $\mathfrak{g}$ as a linear combination of basis $\mathcal B_k$ of $\mathfrak{h}_k$, where $k$ is the output of the routine {\tt DimRepresentation}. Therefore, we start loading Lie algebra $\mathfrak{h}_k$ with the sentence

{\footnotesize
\begin{verbatim}
g > DGsetup(law_h(k),[x],[a])
\end{verbatim}}

In this sentence, {\tt [x]} is used to denote the basis vectors in $\mathcal B_k$ as $\{x_i\}_{i=1}^{\frac{k (k+1)}{2}}$ instead of $\{e_i\}_{i=1}^{\frac{k (k+1)}{2}}$, since we need different notations for the bases of both Lie algebras $\mathfrak{g}$ and $\mathfrak{h}_k$. The notation {\tt [a]} corresponds to the list of structure constants defining the law of Lie algebra $\mathfrak h_n$ and expressed as in (\ref{EcBrack1}), (\ref{EcBrack2}) and (\ref{EcBrack3}).
From this point, we can also work over the Lie algebra $\mathfrak{h}_k$, named {\tt L1} by the package. Next, as we did with Lie algebra $\mathfrak{g}$, we define all the ideals of the derived series of $\mathfrak{h}_k$ as follows


{\footnotesize
\begin{verbatim}
L1 > C[1](h):=[seq(x||i, i=1..k*(k+1)/2)];
L1 > C[2](h):=DerivedAlgebra();
L1 > for i from 3 to nops(DerivedSeries_h(k)) do C[i](h):=DerivedAlgebra(C[i-1](h)); end do;
\end{verbatim}}

%

Next, we express all the vectors $\{e_i\}_{i=1}^n$ from the basis of $\mathfrak{g}$ as a linear combination of basis $\mathcal B_k=\{x_i\}_{i=1}^{\frac{k (k+1)}{2}}$ of $\mathfrak{h}_k$. To do so, we first implement a subroutine called {\tt listposi}, which computes the position of an element within a list. Then, we also implement the routine {\tt expr}. The latter uses the derived series of both $\mathfrak{g}$ and $\mathfrak{h}_k$ and, by applying Proposition~\ref{PropCon}, returns as output two lists: the first contains all the expressions according to the second step of the method, that is, Equation (\ref{comblineal}); and in the second, the conditions over the coefficients so that non-zero vectors are considered.

{\footnotesize
\begin{verbatim}
L1 > listposi:=proc(a,L)
L1 > for i from 1 to nops(L) do
L1 >   if a=L[i] then
L1 >     return i
L1 >   end if;
L1 > end do;
L1 > end proc:
\end{verbatim}}

Let us note that in order to implement routine {\tt expr} we have defined several clusters comparing the derived series of $\mathfrak{g}$ and $\mathfrak{h}_k$. With this rouitne, we express the basis vectors of $\mathfrak{g}$ as a linear combination of the basis vectors from $\mathfrak{h}_k$.


{\footnotesize
\begin{verbatim}
L1 > expr:=proc()
L1 > L:=[];M:=[];
L1 > for i from 1 to index_g - 1 do
L1 >   for j from 1 to nops(C[i](g)) do
L1 >     if member(C[i](g)[j],C[i+1](g))=false then
L1 >       if C[i+1](g)<>[] then
L1 >         L:=[op(L),C[i](g)[j]=sum(a[j,k]*C[i](h)[k],k=1..nops(C[i](h)))];
L1 >         M:=[op(M),sum(a[j,k]^2,k=1..nops(C[i](h)))<>0];
L1 >       end if;
L1 >       else member(C[i](g)[j],C[1](g),'p');
L1 >         N:=[seq(listposi(C[i+1](h)[k],C[1](h)),k=1..nops(C[i+1](h)))];
L1 >         L:=[op(L),C[i](g)[j]=sum(a[p,N[k]]*C[i+1](h)[k],k=1..nops(C[i+1](h)))];
L1 >         M:=[op(M),sum(a[p,N[k]]^2,k=1..nops(C[i+1](h)))<>0];
L1 >     end if;
L1 >   end do;
L1 > end do;
L1 > return L,M;
L1 > end proc:
\end{verbatim}}

After expressing the basis of $\mathfrak{g}$ with respect to basis $\mathcal B_k$ of $\mathfrak{h}_k$, we impose the law of both Lie algebras to the previous expressions. In this way, we implement the routine {\tt Listeq}, which returns two lists: the first one contaning the elements to be equal to zero; and the second one with the conditions to assure the linear independence of the basis.

{\footnotesize
\begin{verbatim}
L1 > Listeq:=proc()
L1 > R:=[];
L1 > for i from 1 to nops(expr()[1])-1 do
L1 >   for j from i+1 to nops(expr()[1]) do
L1 >     if BracketOfSubspaces([lhs(expr()[1][i])],[lhs(expr()[1][j])])=[] then
L1 >       R:=[op(R),op(BracketOfSubspaces([rhs(expr()[1][i])],[rhs(expr()[1][j])]))];
L1 >       else R:=[op(R),op(BracketOfSubspaces([lhs(expr()[1][i])],[lhs(expr()[1][j])]))-
              op(BracketOfSubspaces([rhs(expr()[1][i])],[rhs(expr()[1][j])]))];
L1 >     end if;
L1 >   end do;
L1 > end do;
L1 > return [op(eval(R,expr()[1]))],[op(expr()[2])];
L1 > end proc:
\end{verbatim}}

Next, we define two variables, {\tt Listexp} and {\tt Listcond}, to save the two outputs of the routine {\tt Listeq} respectively. Finally, we implement the routine {\tt sys} to solve the system of equations resulting from the previous expressions. This routine must receive as input the lists {\tt Listexp} and {\tt Listcond}, returning as output the set of solutions which determine the coefficients of the representation of $\mathfrak{g}$ by using Lie algebra $\mathfrak{h}_k$.

{\footnotesize
\begin{verbatim}
L1 > sys:=proc(L,M)
L1 > local Q;
L1 > Q:=[];
L1 > for i from 1 to nops(L) do
L1 >   Q:=[op(Q),seq(coeff(L[i],x||j),j=1..nops(C[1](h)))];
L1 > end do;
L1 > Q:=[op(Q),op(M)];
L1 > return solve(Q);
L1 > end proc:
\end{verbatim}}

If no solution is obtained, then it is not possible to represent $\mathfrak{g}$ as a Lie subalgebra of Lie algebra $\mathfrak{h}_k$ and we must try with the next Lie algebra: $\mathfrak{h}_{k+1}$. Therefore, we would have to repeat the process from the execution of $\mathfrak h_{k+1}$ with the sentence {\tt DGsetup(law$\_$h(k),[x],[a])}, but replacing {\tt k} with {\tt k+1}.
%
%
%

\subsection{Examples of application}

Next, we show an example with the $3$-dimensional solvable Lie algebra with law $[e_1,e_3] = e_2$.
We must run all the routines implemented in the previous section. Here, we only reproduce the most important outputs and those sentences to be modified for this specific example. To define the solvable Lie algebra $\mathfrak{g}$, the follow sentence is run

{\footnotesize
\begin{verbatim}
> L:= _DG([["LieAlgebra", g, [3]], [[[1, 3, 2], 1]]]);
> DGsetup(L);
\end{verbatim}}

We must fill in {\tt DerivedSeries} the value $n=3$ and then, execute the following sentences


{\footnotesize
\begin{verbatim}
> DerivedSeries();
            [3, 1, 0]
> assign(index_g,nops(DerivedSeries()));
> index_g;
            3
> C[1](g):=[seq(e||i, i=1..3)];
            C[1](g):=[e1,e2,e3]
> C[2](g):=DerivedAlgebra();
            C[2](g):=[e2]
> for i from 3 to index_g do C[i](g):=DerivedAlgebra(C[i-1](g));
            C[3](g):=[]
\end{verbatim}}

At this point, loading the routine {\tt DimRepresentation}, we obtain

{\footnotesize
\begin{verbatim}
> DimRepresentation();
          2
\end{verbatim}}

Therefore, we must use $k=2$ in Step 1 of the algorithm. We look for a representation of $\mathfrak{g}$ as a Lie subalgebra of $\mathfrak{h}_2$. Now we start loading Lie algebra $\mathfrak{h}_2$ and we define all the terms of its derived series as follows

{\footnotesize
\begin{verbatim}
g > DGsetup(law_h(2),[x],[a])
          Lie algebra: L1
L1 > MultiplicationTable("LieBracket");
       [[x1, x2] = x2, [x2, x3] = x2]
L1 > C[1](h):=[seq(x||i, i=1..2*3/2)];
           C[1](h):=[x1,x2,x3]
L1 > C[2](h):=DerivedAlgebra();
              C[2](h):=[x2]
L1 > for i from 3 to nops(DerivedSeries_h(2)) do C[i](h):=DerivedAlgebra(C[i-1](h));end do;
              C[3](h):=[]
\end{verbatim}}

By executing the routine {\tt expr}, we obtain the following output

{\footnotesize
\begin{verbatim}
L1 > expr();
     [e1=a[1,1]*x1+a[1,2]*x2+a[1,3]*x3,e2=a[2,2]*x2,e3=a[3,1]*x1+a[3,2]*x2+a[3,3]*x3],
     [a[1,1]^2+a[1,2]^2+a[1,3]^2<>0,a[2,2]^2<>0,a[3,1]^2+a[3,2]^2+a[3,3]^2<>0]
\end{verbatim}}

The output of {\tt Listeq} is

{\footnotesize
\begin{verbatim}
L1 > Listeq();
   [(-a[1,3]*a[2,2]+a[1,1]*a[2,2])*x2,a[2,2]*x2-(-a[1,3]*a[3,2]+a[1,2]*a[3,3]-
   a[1,2]*a[3,1]+a[1,1]*a[3,2])*x2,(a[2,2]*a[3,3]-a[2,2]*a[3,1])*x2],
   [a[1,1]^2+a[1,2]^2+a[1,3]^2<>0,a[2,2]^2<>0,a[3,1]^2+a[3,2]^2+a[3,3]^2<>0]
\end{verbatim}}

After defining the variables {\tt Listexp} and {\tt Listcond} from the previous output, the routine {\tt sys} is executed as follows

{\footnotesize
\begin{verbatim}
L1 > sys(Listexp,Listcond);
\end{verbatim}}

Since no answer is returned, there is no solution for the underlying system. Hence, Lie algebra $\mathfrak{g}$ cannot
be represented as a Lie subalgebra of $\mathfrak{h}_2$. Thus, the process must now be repeated from the execution of the sentence
{\tt DGsetup(law$\_$h(2),[x],[a])}, where $\mathfrak{h}_2$ is replaced with $\mathfrak{h}_3$.

{\footnotesize
\begin{verbatim}
g > DGsetup(law_h(3),[x],[a])
          Lie algebra: L1
L1 > MultiplicationTable("LieBracket");
      [[x1,x2]=x2,[x1,x3]=x3,[x2,x4]=x2,[x2,x5]=x3,[x3,x6]=x3,[x4,x5]=x5,[x5,x6]=x5]
L1 > C[1](h):=[seq(x||i, i=1..3*4/2)];
           C[1](h):=[x1,x2,x3,x4,x5,x6]
L1 > C[2](h):=DerivedAlgebra();
              C[2](h):=[x2,x3,x5]
L1 > for i from 3 to nops(DerivedSeries_h(3)) do
L1 >    C[i](h):=DerivedAlgebra(C[i-1](h));
L1 > end do;
              C[3](h):=[x3]
              C[4](h):=[]
\end{verbatim}}

This time, the variables {\tt Listexp} and {\tt Listcond} are defined
from the output of {\tt Listeq} as follows

{\footnotesize
\begin{verbatim}
L1 > Listexp:=[(-a[1,4]*a[2,2]+a[1,1]*a[2,2])*x2+(-a[1,6]*a[2,3]-a[1,5]*a[2,2]+
     a[1,2]*a[2,5]+a[1,1]*a[2,3])*x3+(-a[1,6]*a[2,5]+a[1,4]*a[2,5])*x5,
     a[2,2]*x2+a[2,3]*x3+a[2,5]*x5-(-a[1,4]*a[3,2]+a[1,2]*a[3,4]-a[1,2]*a[3,1]+
     a[1,1]*a[3,2])*x2-(-a[1,6]*a[3,3]-a[1,5]*a[3,2]+a[1,3]*a[3,6]-a[1,3]*a[3,1]+
     a[1,2]*a[3,5]+a[1,1]*a[3,3])*x3-(-a[1,6]*a[3,5]+a[1,5]*a[3,6]-a[1,5]*a[3,4]+
     a[1,4]*a[3,5])*x5,(a[2,2]*a[3,4]-a[2,2]*a[3,1])*x2+(-a[2,5]*a[3,2]+
     a[2,3]*a[3,6]-a[2,3]*a[3,1]+a[2,2]*a[3,5])*x3+(a[2,5]*a[3,6]-a[2,5]*a[3,4])*x5]:
\end{verbatim}}

{\footnotesize
\begin{verbatim}
L1 > Listcond:=[a[1,1]^2+a[1,2]^2+a[1,3]^2+a[1,4]^2+a[1,5]^2+a[1,6]^2<>0,a[2,2]^2+
     a[2,3]^2+a[2,5]^2<>0,a[3,1]^2+a[3,2]^2+a[3,3]^2+a[3,4]^2+a[3,5]^2+a[3,6]^2<>0];
\end{verbatim}}

Finally, the routine {\tt sys} is executed with {\tt Listexp} and {\tt Listcond} as parameters

{\footnotesize
\begin{verbatim}
L1 > sys(Listexp,Listcond);
{a[1,1]=a[1,6],a[1,2]=a[1,2],a[1,3]=a[1,3],a[1,4]=a[1,6],a[1,5]=a[1,5],a[1,6]=a[1,6],
a[2,2]=0,a[2,3]=-a[1,5]*a[3,2]+a[1,2]*a[3,5],a[2,5]=0,a[3,1]=a[3,1],a[3,2]=a[3,2],
a[3,3]=a[3,3],a[3,4]=a[3,1],a[3,5]=a[3,5],a[3,6]=a[3,1]},
{a[1,1]=a[1,1],a[1,2]=a[1,2],a[1,3]=a[1,3],a[1,4]=a[1,1],a[1,5]=a[1,5],a[1,6]=a[1,1],
a[2,2]=0,a[2,3]=-a[1,5]*a[3,2]+a[1,2]*a[3,5],a[2,5]=0,a[3,1]=a[3,1],a[3,2]=a[3,2],
a[3,3]=a[3,3],a[3,4]=a[3,1],a[3,5]=a[3,5],a[3,6]=a[3,1]}
\end{verbatim}}

As a particular solution of this system, we obtain the representative
$$e_1=-x_5, \quad e_2=x_3, \quad e_3=x_2$$
or, by considering the original notation for basis $\mathcal{B}_3$,
$$e_1 = -X_{2,3}, \quad e_2 = X_{1,3}, \quad e_3 = X_{1,2}$$

To conclude this section, we would like to point out that our algorithmic method is not only valid for low-dimensional solvable Lie algebras; but it provides minimal faithful representations for solvable Lie algebras of higher dimension. In this sense, we have computed some other examples corresponding to solvable Lie algebras such that their dimension is greater than the last being classified, namely: dimension $n\ge 7$. More concretely, Kobel~\cite{Kobel} gave the list of 7-dimensional solvable Lie algebras with codimension 1. From that list, we have considered that with law
$$[e_2,e_7]=e_3, [e_3,e_7]=e_4, [e_4,e_7]=e_5, [e_5,e_7]=e_6, [e_6,e_7]=e_6$$
and our algorithm returns the following representation
$$e_1=X_{1,1}, \, e_2= X_{2,3},\, e_3=X_{2,4},\, e_4=X_{2,5},\, e_5=X_{2,6},$$
$$ e_6=X_{2,6}, \, e_7=X_{3,4}+X_{4,5}+X_{5,6}+X_{6,6}$$

We have also considered a second example consisting of the $8$-dimensional solvable Lie algebra with a 4-dimensional abelian ideal and law
$$[e_1,e_2]=e_3, [e_2,e_5]=e_6, [e_4,e_5]=e_8, [e_1,e_6]=e_7, [e_2,e_6]=e_8, [e_3,e_5]=e_7.$$
In this case, our algorithmic method gives this representation
$$e_1=X_{1,4}+X_{3,5}, \, e_2= X_{1,4}+X_{5,6},\, e_3=X_{3,6},\, e_4=-2X_{1,2},$$
$$e_5=X_{1,3}+X_{2,6}+X_{4,5}, \, e_6=X_{1,5}-X_{4,6}, \, e_7=-X_{1,6}, \, e_8=-2X_{1,6}.$$


\section{Solvable Lie algebras of dimension less than $6$}
\label{Representation}

This section is devoted to apply the algorithm implemented in Section \ref{Algorithmic method}, obtaining a minimal faithful upper-triangular matrix representation for each solvable Lie algebra of dimension less than 6. In addition, we compute such representations for several important families of $n$-dimensional solvable Lie algebras. Tables \ref{class} to \ref{TableClas5DNDRSNNLA2} show the classification of solvable Lie algebras of dimension less than $6$ given in \cite{AND,Mub}, taking into account that we have only written the nonzero brackets in the law of  each Lie algebra; whereas Tables \ref{repres} to \ref{TableRepres5DNDRSNNLA2} contain a representative for each algebra in the previous tables.
In virtue of these tables, we can state the following

\begin{prp}\label{ThrMain1}\quad
A minimal faithful representation by upper-triangular matrices for each solvable Lie algebra of dimension less than $6$ with the dimension of such a minimal representation is given in Tables \ref{repres} to \ref{TableRepres5DNDRSNNLA2}. Moreover, such representations can be obtained with a natural re\-pre\-sen\-ta\-ti\-ve.
\end{prp}

Next, we show several results to determine a representative for minimal faithful upper-triangular matrix representations of three different families of solvable Lie algebras. Proposition \ref{solvablenon-nilpotent} deals with a family of solvable non-nilpotent Lie algebras. Then, Proposition \ref{heisenberg} computes a minimal representative for Heisenberg algebras. These Lie algebras constitute a special subclass of nilpotent Lie
algebras and are very interesting for their applications to both the theory of nilpotent Lie algebras itself and Theoretical Physics. Finally, Proposition \ref{modelfiliform} provides a minimal representation for model filiform Lie algebras. These algebras are the most structured Lie
algebras in the nilpotent class and were introduced by Vergne \cite{Vergne} in 1966.
These propositions can be proved by applying the algorithm considered in Section \ref{Algorithmic method} from a theoretical point of view and not by running the implementation.

\begin{prp}\label{solvablenon-nilpotent}

Let $\mathfrak{s}_n$ be an $n$-dimensional solvable Lie algebra
with basis $\{e_i\}_{i=1}^n$ and law $[e_i,e_n]=e_i$, for $1 \leq i < n$. Then, $\overline{\mu}(\mathfrak{s}_n) = n$. In fact, a
natural representative of $\mathfrak{s}_n$ is given by
$$\{e_{j}=X_{1,j+1}\}_ {j=1}^{j=n-1} \cup \{e_{n}=-X_{1,1}\}$$

\end{prp}

\begin{prp}\label{heisenberg}

Let $\mathfrak{H}_{2n+1}$ be the $(2n+1)$-dimensional Heisenberg algebra with basis $\{e_i\}_{i=1}^{2n+1}$ and law $[e_{2i},e_{2i+1}]=e_1$, for $1 \leq i \leq n$. Then, $\overline{\mu}(\mathfrak{H}_{2n+1}) = n+2$. Moreover, a
natural representative of $\mathfrak{H}_{2n+1}$ is given by

$$\{e_{2j+1}=X_{j+1,n+2}\}_ {j=0}^{j=n} \cup \{e_{2k}=X_{1,k+1}\}_{k=1}^{n+1}$$

%

\end{prp}

\begin{prp}\label{modelfiliform}

Let $\mathfrak{f}_{n}$ be the $n$-dimensional filiform Lie algebra
with basis $\{e_i\}_{i=1}^{n}$ and law $[e_1,e_h]=e_{h-1}$, for $3 \leq h \leq n$. Then, $\overline{\mu}(\mathfrak{f}_{n}) = n$. Moreover, a
natural representative of $\mathfrak{f}_{n}$ is given by
$$\left\{e_1 =\sum_{i=1}^{n-2} X_{i,i+1}\right\}\cup \{e_j = X_{j-1,n}\}_{j=2}^{n}$$

\end{prp}

%

\begin{table}[htp]
\caption{Solvable Lie algebras of dimension less than $5$.}\label{class}
\small
\begin{center}
\begin{tabular}{|c|c|c|}
\hline Dim. & Lie algebra & (Non-zero) Lie brackets \\
\hline
1 & $\mathfrak{s}_1^1$ & --- \\
\hline
2 & $\mathfrak{s}_2^1$ & --- \\
  & $\mathfrak{s}_2^2$ & $[e_1, e_2]\!\! =\! e_1$ \\
\hline
3 & $\mathfrak{s}_3^1$ & --- \\
  & $\mathfrak{s}_3^2$ & $[e_1,e_3]\!\!=\!e_2$ \\
  & $\mathfrak{s}_3^3$ & $[e_1,e_3]\!\! =\! e_1$, \ $[e_2,e_3]\!\! =\! e_2$ \\
  & $\mathfrak{s}_3^4$ & $[e_1,e_3]\!\! =\! e_2$, \ $[e_2,e_3]\!\! =\! -e_1$ \\
  & $\mathfrak{s}_3^5$ & $[e_1,e_3]\!\! =\! -e_1$, \ $[e_2,e_3]\!\!=\! -e_1\!-\!e_2$ \\
  & $\mathfrak{s}_3^6$ & $[e_1,e_3]\!\! =\! -e_1$ \\
\hline
4 & $\mathfrak{s}_4^1$ & --- \\
  & $\mathfrak{s}_4^2$ & $[e_1,e_3]\!\!=\!e_2$ , \ $[e_1,e_4]\!\! =\! e_3$\\
  & $\mathfrak{s}_4^3$ & $[e_1,e_3]\!\! =\! e_3$, \ $[e_1,e_4]\!\! =\! e_4$, \ $[e_2,e_3]\!\! =\! e_4$ \\
  & $\mathfrak{s}_4^4$ & $[e_1,e_3]\!\! =\! e_3$, \ $[e_1,e_4]\!\! =\! e_4$, \ $[e_2,e_3]\!\! =\! -e_4$, \ $[e_2,e_4]\!\! =\! e_3$ \\
  & $\mathfrak{s}_4^5$ & $[e_1,e_3]\!\!=\!e_3$ , \ $[e_1,e_2]\!\! =\! e_4$ \\
  & $\mathfrak{s}_4^6$ & $[e_4,e_1]\!\! =\! e_1$, \ $[e_4,e_2]\!\! =\!\alpha e_2$, \ $[e_4,e_3]\!\! =\!\beta e_3$ \\
  & $\mathfrak{s}_4^7$ & $[e_3,e_1]\!\! =\!\alpha e_1$, \ $[e_3,e_2]\!\! =\! e_2$, \ $[e_3,e_4]\!\! =\!e_2+e_4$   \\
  & $\mathfrak{s}_4^8$ & $[e_1,e_2]\!\! =\! e_2+e_3$, \ $[e_1,e_3]\!\! =\! e_3+e_4$, \ $[e_1,e_4]\!\! =\!e_4$ \\
  & $\mathfrak{s}_4^9$ & $[e_1,e_2]\!\! =\! \beta e_2-e_3$, \ $[e_1,e_3]\!\! =\!e_2+ \beta e_3$, \ $[e_1,e_4]\!\! =\! \alpha e_4$ \\
  & $\mathfrak{s}_4^{10}$ & $[e_1,e_2]\!\! =\! (\alpha-1) e_2$, \ $[e_1,e_3]\!\! =\!e_3$, \ $[e_1,e_4]\!\! =\!\alpha e_4$, \ $[e_2,e_3]\!\! =\! e_4$ \\
  & $\mathfrak{s}_4^{11}$ & $[e_1,e_2]\!\! =\! e_2+e_3$, \ $[e_1,e_3]\!\! =\!e_3$, \ $[e_1,e_4]\!\! =\!2 e_4$, \ $[e_2,e_3]\!\! =\! e_4$\\
  & $\mathfrak{s}_4^{12}$ & $[e_1,e_2]\!\! =\! \alpha e_2-e_3$, \ $[e_1,e_3]\!\! =\!e_2+\alpha e_3$, \ $[e_1,e_4]\!\! =\!2 \alpha e_4$, \ $[e_2,e_3]\!\! =\! e_4$ \\
\hline
\end{tabular}
\end{center}
\end{table}

\begin{table}
\caption{$5$-dimensional non-decomposable real solvable
Lie algebras}\label{TableClas5DNDRSNNLA}
\begin{center}\scriptsize
\begin{tabular}{|c|c|c|}
\hline Lie algebra & (Non-zero) Lie brackets & Parameters  \\
\hline$\mathfrak{g}_{5,1}$ &  $[e_1,e_3]=e_5,\, [e_2,e_4]=e_5$ &   \\
\hline$\mathfrak{g}_{5,2}$ &  $[e_1,e_2]=e_4,\, [e_1,e_3]=e_5$ &  \\
\hline$\mathfrak{g}_{5,3}$ &  $[e_1,e_2]=e_4,\, [e_1,e_4]=e_5,\, [e_2,e_3]=e_5$ &  \\
\hline$\mathfrak{g}_{5,4}$ &  $[e_1,e_2]=e_3,\, [e_1,e_3]=e_4,\, [e_2,e_3]=e_5$ &  \\
\hline$\mathfrak{g}_{5,5}$ &  $[e_1,e_2]=e_3,\, [e_1,e_3]=e_4,\, [e_1,e_4]=e_5$ &  \\
\hline$\mathfrak{g}_{5,6}$ &  $[e_1,e_2]=e_3,\, [e_1,e_3]=e_4,\, [e_1,e_4]=e_5,\, [e_2,e_3]=e_5$ &  \\
\hline $\mathfrak{g}_{5,7}$ & \begin{tabular}{c}$[e_1,e_5]= e_1,
[e_2,e_5]=\alpha e_2,$ \\ $[e_3,e_5]=\beta e_3, [e_4,e_5] = \gamma
e_4$ \end{tabular} &  \begin{tabular}{c} $-1 \leq \gamma \leq \beta \leq \alpha \leq 1$, \\ $\alpha \beta \gamma \neq 0$. \end{tabular} \\
\hline$\mathfrak{g}_{5,8}$ & $[e_2,e_5]= e_1, [e_3,e_5]=e_3, [e_4,e_5]=\gamma e_4,$  &  $0 < | \gamma | \leq 1$  \\
\hline$\mathfrak{g}_{5,9}$ & \begin{tabular}{c}$[e_1,e_5]=e_1, [e_2,e_5]=e_1+e_3, [e_3,e_5]=\beta e_3, [e_4,e_5]=\gamma e_4$\end{tabular}  & $0 \neq \gamma \leq \beta$ \\
\hline$\mathfrak{g}_{5,10}$ & $[e_2,e_5]=e_1, [e_3,e_5]=e_2, [e_4,e_5]=e_4$ &   \\
\hline$\mathfrak{g}_{5,11}$ & \begin{tabular}{c}$[e_1,e_5]=e_1, [e_2,e_5]=e_1+e_2, [e_3,e_5]=e_2 + e_3, [e_4,e_5]=\gamma e_4$\end{tabular} &  $\gamma \neq 0$ \\
\hline$\mathfrak{g}_{5,12}$ & \begin{tabular}{c}$[e_1,e_5]=e_1, [e_2,e_5]=e_1+e_2,$ \\ $[e_3,e_5]=e_2 + e_3, [e_4,e_5]=e_3 + e_4$\end{tabular} &   \\
\hline$\mathfrak{g}_{5,13}$ & \begin{tabular}{c}$[e_1,e_5]=e_1, [e_2,e_5]=\gamma e_2$, \\ $[e_3,e_5]=p e_3 - s e_4, [e_4,e_5]=s e_3 + p e_4$\end{tabular} &  $\gamma s \neq 0, |\gamma|\leq 1$ \\
\hline$\mathfrak{g}_{5,14}$ & \begin{tabular}{c}$[e_2,e_5]=e_1, [e_3,e_5]=p e_3 - e_4, [e_4,e_5]=e_3 + p e_4$ \end{tabular} &   \\
\hline$\mathfrak{g}_{5,15}$ & \begin{tabular}{c}$[e_1,e_5]=e_1, [e_3,e_5]=\gamma e_3,$ \\ $ [e_2,e_5]=e_1+e_2, [e_4,e_5]=e_3 + \gamma e_4$\end{tabular} &  $-1\leq \gamma \leq 1$ \\
\hline$\mathfrak{g}_{5,16}$ & \begin{tabular}{c}$[e_1,e_5]=e_1, [e_2,e_5]=e_1+e_2,$ \\ $[e_3,e_5]=p e_3 - s e_4, [e_4,e_5]=s e_3 + p e_4$\end{tabular} &  $s \neq 0$ \\
\hline$\mathfrak{g}_{5,17}$ & \begin{tabular}{c}$[e_1,e_5]=p e_1 - e_2, [e_2,e_5]=e_1+ p e_2,$ \\ $[e_3,e_5]=q e_3 - s e_4, [e_4,e_5]=s e_3 + q e_4$\end{tabular} &  $s \neq 0$ \\
\hline$\mathfrak{g}_{5,18}$ & \begin{tabular}{c}$[e_3,e_5]=e_1 + p e_3 - e_4, [e_2,e_5]=e_1+ p e_2$ \\ $[e_1,e_5]=p e_1 - e_2, [e_4,e_5]=e_2 + e_3 - p e_4$\end{tabular} &  $p \geq 0$ \\
\hline$\mathfrak{g}_{5,19}$ & \begin{tabular}{c}$[e_2,e_3]=e_1, [e_1,e_5]=(1 + \alpha) e_1,$ \\ $[e_2,e_5]=e_2, [e_3,e_5]=\alpha e_3, [e_4,e_5] = \beta e_4$\end{tabular} &  $\beta \neq 0$ \\
\hline$\mathfrak{g}_{5,20}$ & \begin{tabular}{c}$[e_2,e_3]=e_1, [e_1,e_5]=(1 + \alpha) e_2, [e_2,e_5]=e_2,$ \\ $[e_3,e_5]=\alpha e_3, [e_4,e_5] =e_1 + (1 + \alpha) e_4$\end{tabular} &   \\
\hline$\mathfrak{g}_{5,21}$ & \begin{tabular}{c}$[e_2,e_3]=e_1, [e_1,e_5]=2 e_1, [e_4,e_5] = e_4,$ \\ $[e_2,e_5]=e_2 + e_3, [e_3,e_5]= e_3 + e_4$\end{tabular} & \\
\hline$\mathfrak{g}_{5,22}$ & $[e_2,e_3]=e_1, [e_2,e_5]=e_3, [e_4,e_5] =e_4$  &   \\
\hline$\mathfrak{g}_{5,23}$ & \begin{tabular}{c}$[e_2,e_3]=e_1, [e_1,e_5]=2 e_1, [e_3,e_5]= e_3,$ \\ $[e_2,e_5]=e_2 + e_3, [e_4,e_5] =\beta e_4$\end{tabular} &  $\beta \neq 0$ \\
\hline$\mathfrak{g}_{5,24}$ & \begin{tabular}{c}$[e_2,e_3]=e_1, [e_1,e_5]=2 e_1,[e_3,e_5]= e_3,$ \\ $[e_2,e_5]=e_2 + e_3, [e_4,e_5] =\epsilon e_1 + 2 e_4$\end{tabular} &  $\epsilon = \pm 1$ \\
\hline$\mathfrak{g}_{5,25}$ & \begin{tabular}{c}$[e_2,e_3]=e_1, [e_1,e_5]=2p e_1, [e_4,e_5] =\beta e_4,$ \\ $[e_2,e_5]=pe_2 + e_3, [e_3,e_5]=- e_2 + pe_3$\end{tabular} &  $\beta \neq 0$ \\
\hline$\mathfrak{g}_{5,26}$ & \begin{tabular}{c}$[e_2,e_5]=pe_2 + e_3, [e_1,e_5]=2p e_1, [e_2,e_3]=e_1,$ \\ $ [e_3,e_5]=-e_2 + pe_3, [e_4,e_5] =\epsilon e_1 + 2p e_4$\end{tabular} &  $\epsilon = \pm 1$ \\
\hline$\mathfrak{g}_{5,27}$ & \begin{tabular}{c}$[e_2,e_3]=e_1, [e_3,e_5]= e_3 + e_4,  [e_4,e_5] =e_1 + e_4, [e_1,e_5]=e_1$\end{tabular}  &  \\
\hline$\mathfrak{g}_{5,28}$ & \begin{tabular}{c}$[e_2,e_3]=e_1, [e_2,e_5]=\alpha e_2, [e_4,e_5] = e_4,$ \\ $[e_1,e_5]=(1 + \alpha) e_1, [e_3,e_5]= e_3 + e_4$ \end{tabular}  &  \\
\hline$\mathfrak{g}_{5,29}$ & \begin{tabular}{c} $[e_2,e_3]=e_1, [e_1,e_5]=e_1, [e_2,e_5]=e_2, [e_3,e_5]= e_4$ \end{tabular}  &  \\
\hline$\mathfrak{g}_{5,30}$ & \begin{tabular}{c}$[e_2,e_4]=e_1, [e_3,e_4]=e_2,$ \\ $[e_1,e_5]=(2 + h)e_1, [e_4,e_5]=e_4$, \\ $[e_2,e_5]=(1 + h)e_2, [e_3,e_5]=h e_3$\end{tabular}  &  \\
\hline$\mathfrak{g}_{5,31}$ & \begin{tabular}{c}$[e_2,e_4]=e_1, [e_3,e_4]=e_2, [e_1,e_5]=3e_1, [e_3,e_5]= e_3,$ \\ $[e_2,e_5]=2 e_2, [e_4,e_5]=e_3 + e_4$\end{tabular} &   \\
\hline$\mathfrak{g}_{5,32}$ & \begin{tabular}{c}$[e_2,e_4]=e_1, [e_3,e_4]=e_2,[e_1,e_5]=e_1,$ \\$[e_2,e_5]=e_2,[e_3,e_5]=he_1 + e_3$\end{tabular}  &  \\
\hline$\mathfrak{g}_{5,33}$ & \begin{tabular}{c} $[e_1,e_4]=e_1, [e_3,e_4]=\beta e_3, [e_2,e_5]=e_2, [e_3,e_5]=\gamma e_3$ \end{tabular} & $\beta^2 + \gamma^2 \neq 0$ \\
\hline
\end{tabular}
\end{center}
\end{table}

\begin{table}
\caption{$5$-dimensional non-decomposable real solvable
 Lie algebras (II)}\label{TableClas5DNDRSNNLA2}
\begin{center}\scriptsize
\begin{tabular}{|c|c|c|}
\hline Lie algebra & (Non-zero) Lie brackets &  Parameters \\
\hline$\mathfrak{g}_{5,34}$ & \begin{tabular}{c}$[e_1,e_4]=e_1, [e_2,e_4]=e_2,$ \\ $[e_3,e_4]=e_3,[e_1,e_5]=e_1, [e_3,e_5]= e_2$\end{tabular} & \\
\hline$\mathfrak{g}_{5,35}$ & \begin{tabular}{c}$[e_1,e_4]=h e_1, [e_2,e_4]=e_2,[e_3,e_4]=e_3,$ \\ $[e_2,e_5]=-e_3,[e_1,e_5]=\alpha e_1, [e_3,e_5]= e_2$\end{tabular} &  $h^2 + \alpha^2 \neq 0$ \\
\hline$\mathfrak{g}_{5,36}$ & \begin{tabular}{c}$[e_2,e_3]=e_1, [e_1,e_4]=e_1,$ \\ $[e_2,e_4]=e_2, [e_3,e_5]= e_3,[e_2,e_5]=-e_2$\end{tabular}  &  \\
\hline$\mathfrak{g}_{5,37}$ & \begin{tabular}{c}$[e_2,e_3]=e_1, [e_1,e_4]=2e_1,[e_2,e_4]=e_2,$ \\ $[e_3,e_4]=e_3, [e_2,e_5]=-e_3,[e_3,e_5]=e_2$\end{tabular}  &  \\
\hline$\mathfrak{g}_{5,38}$ &  $[e_1,e_4]=e_1, [e_2,e_5]=e_2,[e_4,e_5] =e_3$ &  \\
\hline$\mathfrak{g}_{5,39}$ & \begin{tabular}{c}$[e_1,e_4]=e_1, [e_2,e_4]=e_2,$ \\ $[e_1,e_5]=-e_2,[e_2,e_5]=e_1,[e_4,e_5]= e_3$\end{tabular}  &  \\
\hline
\end{tabular}
\end{center}
\end{table}

\begin{table}[htp]
\caption{Representation of Solvable Lie algebras of dimension less than $5$.}\label{repres}
\small
\begin{center}

\begin{tabular}{|c|c|c|}
\hline Lie algebra & Representation & $\overline{\mu}$ \\
\hline $\mathfrak{s}_1^1$ & $e_1=X_{11}$ & $1$ \\
\hline $\mathfrak{s}_2^1$ & $e_1=X_{11}, e_2=X_{22}$ & $2$ \\
\hline $\mathfrak{s}_2^2$ & $e_1=X_{11}, e_2=X_{12}$ & $2$\\
\hline $\mathfrak{s}_3^1$ & $e_1=X_{11}, e_2=X_{22}, e_3=X_{33}$ & $3$\\
\hline $\mathfrak{s}_3^2$ & $e_1=X_{12}, e_2 = X_{13}, e_3=X_{23}$ & $3$\\
\hline $\mathfrak{s}_3^3$ & $e_1=X_{13}, e_2=X_{12}, e_3=-X_{1,1}$ & $3$\\
\hline $\mathfrak{s}_3^4$ & $e_1=X_{12}, e_2=i X_{12}+ X_{13}, e_3=iX_{22}+ X_{23}-iX_{33}$ & $3$\\
\hline $\mathfrak{s}_3^5$ & $e_1 = X_{13}, e_2 = X_{23}, e_3 = X_{12} - X_{33}$ & $3$\\
\hline $\mathfrak{s}_3^6$ & $e_1 = X_{12}$, $e_2 = X_{11}+X_{22}$, $e_3= -X_{22}$ & $2$\\
\hline $\mathfrak{s}_4^1$ & $e_1\!\!=\!X_{11}, e_2\!\!=\!X_{22}, e_3\!\!=\!X_{33}, e_4\!\!=\!X_{44}$ & $4$\\
\hline $\mathfrak{s}_4^2$ & $e_1\!\!=\!-(X_{23}+X_{34}), e_2\!\!=\!X_{14}, e_3\!\!=\!X_{13}, e_4\!\!=\!X_{12}$ & $4$\\
\hline $\mathfrak{s}_4^3$ & $e_1\!\!=\!X_{11}, e_2\!\!=\!-X_{23}, e_3\!\!=\!X_{12}, e_4\!\!=\!X_{13}$ & $3$\\
\hline $\mathfrak{s}_4^4$ & $e_1\!\!=\!X_{11}, e_2 \!\!=\! -\!\, i X_{22}\! +\! X_{23}\! +\! i X_{33}, e_3 \!\!=\! X_{12}, e_4 \!\!=\! -\!\, i X_{12}\!+\! X_{13}$ & $3$ \\
\hline $\mathfrak{s}_4^5$ & $e_1\!\!=\! X_{12}\! +\! X_{33}, e_2\!\!=\!X_{24}, e_3\!\!=\!X_{34}, e_4\!\!=\!X_{14}$ & $4$\\
\hline $\mathfrak{s}_4^6$ & $e_1\!\!=\!X_{14}, e_2\!\!=\!X_{24}, e_3\!\!=\!X_{34}, e_4\!\!=\!(\alpha\!-\!1)X_{22}\!+\!(\beta\!-\!1)X_{33}\!-\!X_{44}$ & $4$ \\
\hline $\mathfrak{s}_4^7$ & $e_1\!\!=\!X_{14}$, $e_2\!\!=\!X_{13}$, $e_3\!\!=\!X_{11}\!+\!X_{23}\!+\!(1\!-\!\alpha)X_{44}$, $e_4\!\!=\!-X_{12}$ & $4$\\
\hline $\mathfrak{s}_4^8$ & $e_1\!\!=\!X_{12}\!-\!X_{23}\!-\!X_{44}, e_2\!\!=\!X_{34}, e_3\!\!=\!-X_{24}, e_4\!\!=\!-X_{14}$ & $4$ \\
\hline $\mathfrak{s}_4^9$ & \begin{tabular}{c} $e_1\!\!=\!\alpha X_{11}\!\!+\!(\alpha\!-\!\beta\!\!-\!i)X_{22}\!\!+\!(\beta\!-\!i)X_{33},$ \\ $e_2\!\!=\!i(X_{12}\!\!-\!\!X_{34}), e_3\!\!=\!X_{12}\!\!+\!X_{34}, e_4\!\!=\!X_{14}$ \end{tabular} & $4$ \\
\hline $\mathfrak{s}_4^{10}$ &  $e_1\!\!=\!\alpha X_{11}+X_{22}, e_2\!\!=\!X_{12}, e_3\!\!=\!X_{23}, e_4 \!\!=\! X_{13}$ & $3$ \\
\hline $\mathfrak{s}_4^{11}$ &  {\small $e_1\!\!=\!X_{12}+X_{11}+X_{22}-X_{44}, e_2\!\!=\!X_{23}-X_{34}, e_3\!\!=\!X_{13}, e_4\!\!=\!X_{14}$} & $4$ \\
\hline $\mathfrak{s}_4^{12}$ & \begin{tabular}{c}  $e_1\!\!=\!2 \alpha X_{11}+(\alpha-i)X_{22}$, $e_2\!\!=\!X_{23}+X_{12}$, \\ $e_3\!\!=\!i(X_{23}-X_{12})$, $e_4\!\!=\!2 i X_{13}$ \end{tabular} & $3$ \\
\hline
\end{tabular}
\end{center}
\end{table}

\begin{table}
\caption{Representation of $5$-dimensional non-decomposable real solvable
Lie algebras}\label{TableRepres5DNDRSNNLA}
\begin{center}\scriptsize
\begin{tabular}{|c|c|c|}
\hline Lie algebra & Representation & $\overline{\mu}$ \\
\hline$\mathfrak{g}_{5,1}$ &  $e_1=X_{1,2}, e_2=X_{1,3}, e_3=X_{1,4}, e_4=X_{1,5}, e_6=X_{2,3}+X_{3,4}+X_{4,5}$ & $4$  \\
\hline$\mathfrak{g}_{5,2}$ &  $e_1=X_{1,2}+X_{2,4}+X_{3,5}, e_2=X_{1,3}, e_3=X_{1,4}, e_4=X_{1,5}, e_5=X_{2,3}+X_{3,4}+X_{4,5}$ & $5$ \\
\hline$\mathfrak{g}_{5,3}$ &  $e_1=X_{1,2}-X_{3,5}, e_2=X_{1,3}+X_{2,5}, e_3=X_{1,4}, e_4=X_{1,5}, e_5=X_{2,3}+X_{3,4}$ & $5$ \\
\hline$\mathfrak{g}_{5,4}$ & $e_1=X_{1,2}+X_{2,3}, e_2=X_{1,3}, e_3=X_{1,4}, e_4=X_{2,4}, e_5=X_{3,4}$  & $4$ \\
\hline$\mathfrak{g}_{5,5}$ & $e_1=X_{1,2}, e_2=X_{1,3}, e_3=X_{1,4}, e_4=X_{2,3}, e_5=X_{2,4}$ & $4$ \\
\hline$\mathfrak{g}_{5,6}$ & $e_1=X_{1,2}, e_2=X_{1,3}, e_3=X_{1,4}, e_4=X_{2,4}, e_5=X_{3,4}$ & $4$  \\
\hline $\mathfrak{g}_{5,7}$ & \begin{tabular}{c} $e_1=X_{1,4}, e_2=X_{1,2}, e_3=X_{1,5}, e_4=X_{3,5},$ \\
$e_5=\alpha X_{2,2}+(\beta-\gamma)X_{3,3}+X_{4,4}+\beta X_{5,5}$ \end{tabular} & $5$ \\
\hline$\mathfrak{g}_{5,8}$ & \begin{tabular}{c} $e_1=X_{1,4}, e_2=X_{1,2}+X_{1,4}, e_3=X_{1,5}, e_4=X_{1,3},$ \\
$e_5= X_{2,4}+ \gamma X_{3,3}+X_{5,5}$ \end{tabular} & $5$ \\
\hline$\mathfrak{g}_{5,9}$ & \begin{tabular}{c} $e_1=X_{2,4}, e_2=X_{1,3}+X_{2,4}, e_3=\beta X_{1,3}, e_4=X_{2,3},$ \\
$e_5= (\beta-\gamma)X_{2,2}+ \beta X_{3,3}+(1+ \beta - \gamma)X_{4,4}$ \end{tabular} & $4$ \\
\hline$\mathfrak{g}_{5,10}$ & \begin{tabular}{c} $e_1=-X_{1,5}, e_2=-X_{1,4}, e_3=X_{1,2}-X_{1,3}, e_4=X_{1,3}+X_{1,4}+X_{1,5},$ \\
$e_5=X_{1,1}+X_{2,2}+X_{2,3}+2 X_{3,3}+X_{3,4}+X_{4,4}+X_{4,5}+X_{5,5}$ \end{tabular} & $5$ \\
\hline$\mathfrak{g}_{5,11}$ & \begin{tabular}{c} $e_1=X_{1,5}, e_2=\frac{1}{2}X_{1,4}+X_{2,5}, e_3=X_{1,5}+X_{2,4}, e_4=X_{3,5},$ \\
$e_5=-\frac{1}{2}X_{1,2}+(1-\gamma)X_{3,3}+X_{4,4}+X_{4,5}+X_{5,5}$ \end{tabular} & $5$ \\
\hline$\mathfrak{g}_{5,12}$ & \begin{tabular}{c} $e_1=X_{1,5}, e_2=X_{1,4}, e_3=X_{1,3}, e_4=X_{1,2},$ \\
$e_5=X_{2,2}+X_{2,3}+X_{3,3}+X_{3,4}+X_{4,4}+X_{4,5}+X_{5,5}$ \end{tabular} & $5$ \\
\hline$\mathfrak{g}_{5,13}$ & \begin{tabular}{c} $e_1=X_{2,4}, e_2=X_{3,4}, e_3=X_{1,4}, e_4=i X_{1,4}$ \\ $e_5=(p-1-s i)X_{2,2}+(p-\gamma-s i)X_{3,3}+(p-s i)X_{4,4}$  \end{tabular} & $4$ \\
\hline$\mathfrak{g}_{5,14}$ & \begin{tabular}{c} $e_1=X_{2,4}, e_2=X_{2,3}, e_3=X_{1,4}, e_4=-iX_{1,4},$ \\
$e_5=iX_{2,2}+iX_{3,3}+X_{3,4}+iX_{4,4}$ \end{tabular} & $4$ \\
\hline$\mathfrak{g}_{5,15}$ & $e_1=-X_{1,4}, e_2=X_{2,4}, e_3=X_{1,3}, e_4=-X_{2,3}, e_5=X_{1,2}+\gamma X_{3,3}+X_{4,4}$ & $4$ \\
\hline$\mathfrak{g}_{5,16}$ & \begin{tabular}{c} $e_1=X_{2,4}, e_2=X_{2,3}, e_3=X_{1,4}, e_4=-i X_{1,4}$ \\ $e_5=s i X_{2,2}+(1+s i)X_{3,3}+X_{3,4}+(1+s i)X_{4,4}$  \end{tabular} & $4$ \\
\hline$\mathfrak{g}_{5,17}$ & \begin{tabular}{c} $e_1=X_{2,4}, e_2=i X_{2,4}, e_3=X_{1,4}, e_4=i X_{1,4}$ \\
$e_5=(-p+i q -s i) X_{2,2}+(q-s i)X_{4,4}$  \end{tabular} & $4$ \\
\hline$\mathfrak{g}_{5,18}$ & $e_1=X_{2,4}, e_2=i X_{2,4}, e_3=-\frac{i(i-p)}{p} X_{2,4}, e_4=\frac{i}{p}X_{2,4}, e_5=(p-i)X_{4,4}$ & $4$\\
\hline$\mathfrak{g}_{5,19}$ & \begin{tabular}{c} $e_1=X_{1,4}, e_2=X_{1,2}, e_3=X_{2,4}, e_4=X_{3,4}$ \\ $e_5=X_{2,2}+(1+\alpha-\beta) X_{3,3}+ \sqrt{\beta(2-\beta)+2\alpha(\beta-\alpha-2)-3} X_{3,4}+(1+\alpha)X_{4,4}$ \end{tabular} & $4$ \\
\hline$\mathfrak{g}_{5,20}$ & \begin{tabular}{c} $e_1=-X_{1,4}, e_2=-X_{1,3}, e_3=Z X_{2,4}+X_{3,4}, e_4=X_{2,4}$ \\ $e_5=-(1+\alpha)X_{1,1}+X_{2,2}-Z X_{1,3}-(1+\alpha)X_{2,2}+Z X_{2,3}+X_{2,4}-\alpha X_{3,3}$ \end{tabular}  & $4$ \\
\hline$\mathfrak{g}_{5,21}$ & \begin{tabular}{c} $e_1=X_{1,4}, e_2=X_{2,3}+X_{3,4}, e_3=-X_{1,3}+X_{2,5}, e_4=-2 X_{1,5}$ \\ $e_5=X_{1,2}+X_{3,3}+X_{3,5}+2X_{4,4}+X_{5,5}$ \end{tabular} & $5$ \\
\hline$\mathfrak{g}_{5,22}$ & $e_1=X_{1,5}, e_2=X_{1,2}+X_{2,3}, e_3=X_{2,5}, e_4=X_{1,4}, e_5=X_{3,5}+X_{4,4}$ & $5$ \\
\hline$\mathfrak{g}_{5,23}$ & \begin{tabular}{c} $e_1=X_{1,4}, e_2=X_{2,3}+X_{3,4}, e_3=-X_{1,3}, e_4= X_{1,5}$ \\ $e_5=-X_{1,1}+X_{1,2}-X_{2,2}+X_{4,4}+(\beta-1)X_{5,5}$ \end{tabular} & $5$ \\
\hline$\mathfrak{g}_{5,24}$ &  \begin{tabular}{c} $e_1=X_{1,4}, e_2=X_{2,3}+X_{3,4}, e_3=-X_{1,3}, e_4= X_{1,5}-\epsilon X_{2,4}$ \\ $e_5=-X_{1,1}+X_{1,2}-X_{2,2}+X_{4,4}+X_{5,5}$ \end{tabular} & $5$ \\
\hline$\mathfrak{g}_{5,25}$ & \begin{tabular}{c} $e_1=-2 iX_{1,4}, e_2=-i X_{1,3}+ i X_{3,4}, e_3=X_{1,3}+X_{3,4}, e_4=(\beta-2p)X_{1,2}+Z X_{1,4}$ \\ $e_5=-2pX_{1,1}+(\beta-2p)X_{2,2}+Z X_{2,4}+(i-p)X_{3,3}, Z=\sqrt{1+2pi+4p\beta-\beta^2-9p^2}$  \end{tabular} & $4$ \\
\hline$\mathfrak{g}_{5,26}$ & \begin{tabular}{c} $e_1=2 \omega X_{1,4}, e_2=\omega X_{1,2}-(1-p^2+pi)X_{1,4}-\omega X_{2,4}, e_3=X_{1,2}+(-p+i)X_{1,4}+X_{2,4},$ \\ $e_4=-2\epsilon \omega X_{3,4}, e_5=X_{1,3}+(p-\omega)X_{2,2}+X_{2,4}+ZX_{3,4}+2p X_{4,4},$\\
  $Z=\sqrt{2p^3i+2pi-2p^4-5p^2-1}, \omega=-p^3+i(p^2+1)$   \end{tabular} & $4$ \\
\hline$\mathfrak{g}_{5,27}$ & $e_1=-X_{1,4}, e_2=X_{1,3}, e_3=-X_{3,4}, e_4=X_{2,4}, e_5=X_{1,2}+X_{2,3}+X_{4,4}$ & $4$ \\
\hline$\mathfrak{g}_{5,28}$ & \begin{tabular}{c} $e_1=-X_{1,4}, e_2=(1+\alpha)iX_{1,4}+X_{2,4}, e_3=X_{1,2}, e_4=X_{1,3},$ \\ $e_5=(1+\alpha)iX_{1,2}+X_{2,2}+X_{2,3}+X_{3,3}+(1+\alpha)X_{4,4}$ \end{tabular} & $4$ \\
\hline$\mathfrak{g}_{5,29}$ & \begin{tabular}{c} $e_1=X_{1,4}, e_2=X_{1,3}, e_3=X_{1,1}+X_{2,2}+X_{3,3}+X_{3,4}+X_{4,4},$ \\ $e_4=-X_{2,4}, e_5=X_{2,2}+X_{2,3}+X_{3,3}+X_{4,4}$ \end{tabular} & $4$\\
\hline$\mathfrak{g}_{5,30}$ & \begin{tabular}{c} $e_1=-2X_{1,4}, e_2=-X_{1,3}+X_{1,4}+X_{2,4}, e_3=X_{2,3}-X_{4,4},$ \\ $e_4=X_{1,2}+X_{3,4}, e_5=X_{2,2}+(1+h)X_{3,3}+X_{3,4}+(2+h)X_{4,4}$ \end{tabular} & $5$ \\
\hline$\mathfrak{g}_{5,31}$ & \begin{tabular}{c} $e_1=\frac{3}{2}X_{1,5}, e_2=X_{1,4}-X_{2,5}, e_3=\frac{1}{3}X_{1,3}-X_{2,4},$ \\ $e_4=\frac{1}{2}X_{1,2}+X_{3,4}+X_{4,5}, e_5=X_{2,2}+X_{2,3}+X_{3,3}+2X_{4,4}+3 X_{5,5}$ \end{tabular} & $5$ \\
\hline
\end{tabular}
\end{center}
\end{table}

\begin{table}
\caption{Representation of $5$-dimensional non-decomposable real solvable
 Lie algebras (II)}\label{TableRepres5DNDRSNNLA2}
\begin{center}\scriptsize
\begin{tabular}{|c|c|c|}
\hline Lie algebra & Representation & $\overline{\mu}$ \\
\hline$\mathfrak{g}_{5,32}$ & \begin{tabular}{c} $e_1=X_{1,4}, e_2=X_{1,3}, e_3=X_{1,2}, e_4=X_{2,2}+X_{3,4},$ \\ $e_5=X_{2,2}+hX_{2,4}+X_{3,3}+X_{4,4}$ \end{tabular}  & $4$\\
\hline$\mathfrak{g}_{5,33}$ & \begin{tabular}{c} $e_1=X_{2,4}, e_2=X_{1,4}, e_3=\gamma X_{2,3}+X_{2,4},$ \\ $e_4=X_{1,1}+\beta X_{3,3} + \frac{\beta-1}{\gamma} X_{3,4} + X_{4,4}, e_5=-X_{1,1}+\gamma X_{3,3}+X_{3,4}$ \end{tabular} & $4$ \\
\hline$\mathfrak{g}_{5,34}$ & $e_1=X_{1,2}, e_2=-X_{1,4}, e_3=X_{3,4}, e_4=X_{2,2}+X_{4,4}, e_5=-X_{1,1}+X_{1,3}-X_{3,3}-X_{4,4}$ & $4$ \\
\hline$\mathfrak{g}_{5,35}$ &  $e_1=X_{2,4}, e_2=X_{1,3}, e_3=-iX_{1,3}, e_4=X_{3,3}+ h X_{4,4}, e_5=-\alpha X_{2,2}+iX_{3,3}$ & $4$ \\
\hline$\mathfrak{g}_{5,36}$ & $e_1=-X_{1,4}, e_2=X_{2,4}, e_3=X_{1,2}, e_4=X_{4,4}, e_5= X_{2,2}+iX_{3,3}$ & $4$ \\
\hline$\mathfrak{g}_{5,37}$ & \begin{tabular}{c} $e_1=2 i X_{1,4}, e_2=X_{1,2}-iX_{1,4}+X_{2,4}, e_3=-iX_{1,2}-X_{1,4}+iX_{2,4},$ \\ $ e_4=X_{2,2}+iX_{2,4}+2X_{4,4}, e_5= iX_{2,2}+X_{2,4}$ \end{tabular} & $4$ \\
\hline$\mathfrak{g}_{5,38}$ & \begin{tabular}{c} $e_1=X_{1,5}, e_2=X_{1,4}, e_3=X_{1,3}, e_4=X_{1,2}+X_{5,5}$ \\
$e_5=X_{1,1}+X_{2,2}+X_{2,3}+X_{3,3}+2 X_{4,4}+X_{5,5}$ \end{tabular} & $5$  \\
\hline$\mathfrak{g}_{5,39}$ & $e_1=X_{2,4}, e_2=-i X_{2,4}, e_3=X_{1,4}, e_4=-X_{2,2}+X_{3,4}, e_5=i X_{1,1}-X_{1,3}+i X_{3,3} + iX_{4,4}$ & $4$  \\
\hline
\end{tabular}
\end{center}
\end{table}


\begin{thebibliography}{99}

\bibitem{Jac} N. Jacobson, A note on automorphisms and
derivations of Lie algebras, {\it Proc. Amer. Math. Soc.} {\bf 6}
(1955) 281--283.

\bibitem{Ner} Y.A. Neretin, A construction of finite-dimensional faithful representation of Lie algebras, {\it Rend. Circ. Mat. Palermo Supp.} {\bf 71} (2003) 159--161.

\bibitem{Var98} V.S. Varadarajan, {\it Lie Groups, Lie Algebras and
their Representations}, Selected Monographies {\bf 17}, Coll\ae ge
Press, Beijing, 1998.

\bibitem{Bu98} D. Burde, On a refinement of Ado's Theorem. {\it
Arch. Math. (Basel)} {\bf 70} (1998) 118--127.

\bibitem{G} R. Ghanam, I. Strugar, G. Thompson, Matrix representations for low dimensional Lie algebras, {\it Extracta Math.} {\bf 20} (2005) 151--184.

\bibitem{Mil77} J. Milnor, On fundamental groups of complete affinely flat manifolds, {\it Adv. Math.} {\bf 25} (1977) 178--187.

\bibitem{Aus64} L. Auslander, The structure of complete locally affine manifolds, {\it Topology} {\bf 3}: Suppl. 1 (1964) 131--139.

\bibitem{Aus77} L. Auslander, Simply transitive groups of affine motions, {\it Amer. J. Math.} {\bf 99}:4 (1977) 809--826.

\bibitem{B} J.C. Benjumea, F.J. Echarte, J. N\'{u}\~{n}ez, A.F. Tenorio, A method to obtain the Lie group associated with a
nilpotent Lie algebra, {\it Comput. Math. Appl.} {\bf 51} (2006) 1493--1506.


\bibitem{MathScan08} J.C. Benjumea, J. N\'{u}\~{n}ez, A.F. Tenorio, Minimal linear representations of the low-dimensional
nilpotent Lie algebras, {\it Math. Scand.} {\bf 102} (2008) 17--26.

\bibitem{CNT} M. Ceballos, J. N\'{u}\~{n}ez, A.F. Tenorio, Representing Filiform Lie Algebras Minimally and Faithfully by Strictly Upper-Triangular Matrices, {\it Journal of Algebra and its Applications} {\bf  12} (2013) 1250196.

\bibitem{NuTeMega07} J. N\'u\~nez, A.F. Tenorio, Minimal Faithful Upper-triangular Matrix Representations for Low-Dimensional Solvable Lie Algebras, in: Electronic Proceedings of MEGA 2007 - Effective Methods in Algebraic Geometry, Strobl, 2007, 13 pp.

\bibitem{AND} A. Andrada, M.L. Barberis, I.G. Dotti, G.P.
Ovando, Pro\-duct structures on four dimensional solvable Lie
algebras, {\it Homology Homotopy Appl.} \textbf{7}
(2005) 9--37.

\bibitem{Mub} G.M. Mubarakzyanov, On solvable Lie algebras, {\it
Izv. Vyssh. Uchebn. Zaved. Mat.} \textbf{32} (1) (1963) 114--123
(in Russian).

\bibitem{Mub1} G.M. Mubarakzyanov, The classification of the real structure
of five-dimensional Lie algebras, {\it Izv. Vyssh. Uchebn. Zaved.
Mat.} \textbf{34} (3) (1963) 99--106 (in Russian).


\bibitem{Turi} P. Turkowski,. Solvable Lie algebras of dimension
six, {\it J. Math. Phys.} \textbf{31} (1990) 1344--1350.

\bibitem{Vergne} M. Vergne, Cohomologie des alg\`{e}bres de Lie nilpotentes, Application \`{a} l'\'{e}tude de la vari\'{e}t\'{e} des algebres de Lie
nilpotentes, {\it Bull. Soc. Math. France} {\bf 98} (1970) 81--116.


\bibitem{F-H91} W. Fulton, J. Harris, {\it Representation theory: a first course}, Springer-Verlag, New York,
1991.

\bibitem{Kobel} C. Kobel, {\it On the Classification of Solvable Lie Algebras of Finite Dimension Containing an Abelian Ideal of
Codimension One}, Master's Thesis, School of Information Science, Computer and Electrical Engineering,
Halmstad University, Halmstad, 2008.


\end{thebibliography}
\end{document}